\documentclass[11pt, a4paper]{article}

\usepackage{amssymb}
\usepackage{amsmath}
\usepackage{amscd}
\usepackage{theorem}
\usepackage{mathrsfs}
\usepackage{enumerate}

\theoremstyle{change}

\newtheorem{thm}{Theorem.}
\newtheorem{lem}[thm]{Lemma.}
\newtheorem{prop}[thm]{Proposition.}

\newtheorem{cor}[thm]{Corollary.}

\newcommand{\T}{\mathbb{T}}
\newcommand{\C}{\mathbb{C}}
\newcommand{\D}{\mathbb{D}}

\newcommand{\B}{\mathbb{B}}

\newcommand{\Hcal}{\mathcal{H}}
\newcommand{\Scal}{\mathcal{S}}

\newcommand{\Tcal}{\mathcal{T}}
\newcommand{\Kcal}{\mathcal{K}}

\newcommand{\Rcal}{\mathcal{R}}

\newcommand{\trace}{\text{tr}}

\renewcommand{\phi}{\varphi}

\newcommand{\proof}{\textbf{Proof. }}
\newcommand{\proofend}{\hfill $\Box$}

\newcommand{\ran}{\mathrm{ran\,}}

\begin{document}
\begin{center}
{\large{\textbf{On Schatten-class perturbations of Toeplitz operators}}}
\end{center}

\vspace{1cm}

\centerline{Michael Didas, J\"org Eschmeier, Dominik Schillo}

\vspace{.5cm}

\begin{center}
\parbox{12cm}{\small Let $D$ denote the unit ball or the unit polydisc in $\C^n$ with $n \geq 2$.  
For $1 \leq p \leq 2n$ in the case of the ball and $1 \leq p < \infty$ for the polydisc, 
we show that a bounded operator $S$ on the Hardy space $H^2(D)$ commutes with all analytic Toeplitz operators modulo
the Schatten class $\Scal_p$ if and only if $S=X+K$ with an analytic Toeplitz operator $X$ and an operator
$K \in \Scal_p$.
This partially answers a question of Guo and Wang \cite{guo-wang}. 
For $1 \leq p < \infty$ and a strictly pseudoconvex or bounded symmetric and circled domain $D \subset \mathbb C^n$, we show that
a given operator $S$ on $H^2(D)$ is a Schatten-$p$-class perturbation of a Toeplitz operator if and only if
$T^*_{\theta} S T_{\theta} - S \in \mathcal S_p$ for every inner function $\theta$ on $D$.
\vspace{0.5cm}

}

\end{center}
\vspace{1cm}

\centerline{\textbf{\S1 \, Introduction and main results}}

\vspace{.5cm}

Let $D \subset \C^n$ be either the unit ball $\B_n$ or the unit polydisc $\D^n$.
In \cite{guo-wang} Guo and Wang asked for a characterization of those bounded operators $S$ on the Hardy space $H^2(D)$
that commute with all analytic Toeplitz operators modulo the Schatten classes
\[
\Scal_p = \{ X \in B(H^2(D)): \|X\|_p = (\trace|X|^p)^{1/p} <\infty \}.
\]
Note that there is a chain of inclusions
\[
\Scal_0 \subset \Scal_p \subset \Scal_q \subset \Scal_\infty \qquad (1\leq p<q<\infty),
\]
where $\Scal_0$ stands for the finite-rank and $\Scal_\infty$ for the compact operators on $H^2(D)$.
Questions concerning the $p$-essential commutant 
\[
\Tcal_a^{\text{ec},p} = \{ S \in B(H^2(D)): [S,X] = SX - XS \in \Scal_p \text{ for every } X \in \Tcal_a\}
\]
of the set of all analytic Toeplitz operators $\Tcal_a = \{ T_f : f \in H^\infty(D) \}$
date back at least to the work of Douglas and Sarason on Toeplitz operators on the unit disc. A complete characterization of the essential commutant
$\Tcal^{ec}_a  = \Tcal^{ec,\infty}_a$ of the analytic Toeplitz operators on the unit disc was given by Davidson \cite{davidson} in 1977. 
A generalization of Davidson's result to the unit ball was obtained in 1997 by Ding and Sun \cite{ding-sun}
(see also Le \cite{trieu-le} for related work). In a recent paper of Everard and the first two authors
\cite{esscom} Davidson's result is extended to the case of smoothly bounded strictly pseudoconvex domains
$D \subset \C^n$. More precisely, it is shown that
\[
\Tcal_a^{\text{ec}} = \{ T_f + K : K \text{ compact},\, f \in L^\infty(\sigma) \text{ with } H_f \text{ compact}  \},
\]
where $\sigma$ denotes the normalized surface measure
on $\partial D$ and $H_f$ is the Hankel operator with symbol $f$.
On the unit disc, a well known theorem of Hartman says that $H_f$ is compact if and only if $f \in H^\infty+C$ and thus
we recover the main result of \cite{davidson}.

At the other end of the chain, namely for $p=0$, Guo and Wang completely solved both the ball and the polydisc case
extending corresponding one-dimensional results of Gu \cite{gu}. They proved in \cite{guo-wang} that
\[
\Tcal_a^{\text{ec},0} = \Tcal_a + \Scal_0 \qquad (\text{if } D = \B_n \text{ or } D = \D^n \text{ and } n\geq 2)
\]
while, for $n=1$, $\Tcal_a$ has to be replaced by the set of all Toeplitz operators $T_f$ where $f$ is the sum of a rational function
and a bounded analytic function.

One of the goals of the present paper is to fill in at least some of the gap between $p=0$ and $p=\infty$.
Using results on the compactness and the Schatten-class membership of Hankel operators on $D = \D^n$ and $D = \B_n$
due to Cotlar and Sadoski \cite{cotlar} and Fang and Xia \cite{fang-xia}, respectively,  we are able prove the following result:

\begin{thm} \label{esscom_p_ball_polydisc}
Let $n \geq 2$ be given. If either
\[
D = \B_n \text{ and } 1\leq p \leq 2n \qquad \text{or} \qquad D = \D^n \text{ and } 1 \leq p < \infty,
\]
then the identity
\[
\Tcal_a^{\text{ec},p} = \Tcal_a + \Scal_p
\]
holds on the Hardy space $H^2(D)$.
\end{thm}


A result of Xia from \cite{xia-characterization}, answering a question of Douglas \cite{douglas}, shows that a given operator
$S \in B(H^2(\mathbb D))$ on the Hardy space of the unit disc is a compact perturbation of a Toeplitz operator 
$T_f$ with symbol $f \in L^{\infty}(\partial \mathbb D)$ if and only if $T^*_{\theta}S T_{\theta} - S$ is compact for every
inner function $\theta$. The question whether a corresponding result holds in higher dimensions seems still to be
open. In the present paper we show that on the Hardy space $H^2(D)$ over the Shilov boundary of a smooth strictly
pseudoconvex or bounded symmetric and circled domain $D \subset \mathbb C^n$ the corresponding characterization
of Schatten-p-class perturbations of Toeplitz operators holds true for $1 \leq p < \infty$.

\begin{thm}\label{xia}
Let $D \subset \mathbb C^n$ be a strictly pseudoconvex domain with smooth boundary or a bounded symmetric
and circled domain and let $1 \leq p < \infty$. Then a given operator $S \in B(H^2(D))$ is of the form
\[
S = T_f + K,
\]
where $T_f $ is a Toeplitz operator with symbol $f \in L^{\infty}(\sigma)$ and $K \in \mathcal S_p(H^2(D))$, if and only if
$T^*_{\theta}S T_{\theta} - S \in \mathcal S_p(H^2(D))$ for every inner function $\theta$ on $D$.
\end{thm}

We prove the above results in the more general framework of regular $A$-isometries. 
A detailed study of these operator tuples and their associated Toeplitz operators can be found in the recent paper \cite{toeplitz-projection}
of Everard and the second author. For the convenience of the reader, we briefly outline the basic facts.
Let $K\subset \C^n$ be a compact set and let $A$ be  a closed subalgebra of the algebra $C(K)$ of all continuous complex-valued functions
on $K$ such that $A \supset \C[z_1, \ldots, z_n]$. We denote by $\partial_A$ the Shilov boundary of $A$ and write
$M^+_1(\partial_A)$ for the set of all regular probability measures on $\partial_A$.
A subnormal tuple $T \in B(\Hcal)^n$ on a complex Hilbert space $\Hcal$ is called an $A$-isometry \cite{esscom} if the spectrum
of the minimal normal extension $U \in B(\Kcal)^n$ of $T$ is contained in $\partial_A$ and if
$A$ is contained in the restriction algebra $\Rcal_T = \{ f \in L^\infty(\mu) : \Psi_U(f)\Hcal \subset \Hcal \}$.
Here $\mu \in M^+_1(\partial_A)$ is a fixed scalar spectral measure of $U$ and $\Psi_U:L^\infty(\mu) \to B(\Kcal)$ denotes
the associated $L^\infty$-functional calculus of $U$.

An $A$-isometry $T \in B(\Hcal)^n$ is called regular if the triple $(A|\partial_A, \partial_A, \mu)$
is regular in the sense of Aleksandrov \cite{aleks} which means by definition that, for every continuous function $\phi \in C(\partial_A)$
with $\phi>0$ on $\partial_A$, there is a sequence of functions $(\phi_k)_{k\geq 1}$ in $A$ such that
 $|\phi_k|<\phi$ $(k \geq 1)$ on $\partial_A$ 
and $\lim_{k \to \infty} |\phi_k| = \phi$ $\mu$-a.e. on $\partial_A$. In his celebrated work \cite{aleks}, Aleksandrov proved that
under this regularity assumption the abstract inner function problem related to $A$ is solvable.
This will be formulated in a more precise way below.
 
It turns out that regular $A$-isometries 
provide a natural and unified framework to treat some
operator-theoretic questions on Hardy spaces over strictly pseudoconvex and bounded symmetric domains $D \subset \C^n$.
More explicitly, taking $A$ to be the domain algebra 
\[
A(D) = \{ f \in C(\overline{D}): f|D \text{ is holomorphic}\}
\]
one can show for both classes of domains that with $A = A(D)$ the triple $(A|\partial_A, \partial_A, \mu)$ is regular
 whenever $\mu$ is a finite positive regular Borel measure on the Shilov boundary $\partial_A$
(see Aleksandrov \cite{aleks}, or \cite{diss} and \cite{subnormal} for more detailed arguments).
For each such measure $\mu$, the abstract Hardy-space multiplication tuple
\[
  M_z = (M_{z_1}, \ldots, M_{z_n}) \in B(H^2_A(\mu))^n \quad \text{with} \quad H^2_A(\mu) = \overline{A|\partial_A}^{\|\cdot\|_{2,\mu}}
	\subset L^2(\mu)
\]
is a regular $A$-isometry. By choosing $\mu$ as the canonical probability measure $\sigma$ on $\partial_{A(D)}$
we obtain the usual Hardy spaces on the unit ball and polydisc.

Theorem 1 is obtained as a corollary of a more general result for regular $A$-isometries.
To formulate this result, we need some more notations: Given a regular $A$-isometry $T\in B(\Hcal)^n$,
let $\Psi_U: L^\infty(\mu) \to B(\Kcal)$ as above denote the functional calculus of the minimal normal extension $U$ of $T$.
Via the formula
\[
T_f = P_{\Hcal} \Psi_U(f) |\Hcal \qquad (f \in L^\infty(\mu))
\]
one can define Toeplitz operators with $L^\infty(\mu)$-symbols associated with $T$. 
Toeplitz operators of this type will be called concrete Toeplitz operators.
The set of analytic $T$-Toeplitz operators is defined as
\[
\Tcal_a(T) = \{ T_f : f \in H^\infty_A(\mu) \}  \subset B(\Hcal) \quad \text{where} \quad H^\infty_A(\mu) = \overline{A}^{w^*} \subset L^\infty(\mu).
\]
Since $A$ is assumed to be regular, Aleksandrov's work \cite{aleks} on the existence of  abstract inner functions guarantees that the set
\[
I_\mu = \{ \theta \in H^\infty_A(\mu) : |\theta| = 1 \quad \mu-\text{a.e. on } \partial_A \}
\]
of all $\mu$-inner functions generates $L^\infty(\mu)$ as a von Neumann algebra (Coroallary 2.5 in  \cite{inner-spherical}).
Let us denote by $\Tcal(T) \subset B(\Hcal)$ the set of all operators $X \in B(\Hcal)$ satisfying
the Brown-Halmos type condition
\[
T_{\overline{\theta}} X T_\theta = X \qquad \text{for every } \theta \in I_\mu .
\]
We shall refer to these operators as abstract $T$-Toeplitz operators. Obviously $\Tcal(T)$ contains each concrete
$T$-Toeplitz operator $T_f$ with symbol $f \in L^{\infty}(\mu)$. 
Using results of Prunaru \cite{prunaru} on spherical isometries one can show that
\[
\Tcal(T) = P_\Hcal (U)'| \Hcal,
\]
where $(U)'$ denotes the commutant of $U$ (see \cite{toeplitz-projection} and the references therein).
With these definitions the announced theorem on $A$-isometries reads as follows.

\begin{thm}\label{esscom_p} Let $T \in B(\Hcal)^n$ be a regular $A$-isometry and $1\leq p < \infty$.
Then every operator $S \in \Tcal_a(T)^{\text{ec},p}$ has the form
\[
S= X + K \quad \text{with} \quad X \in \Tcal(T) \text{ and } K \in \Scal_p(\Hcal). 
\]
\end{thm}
In contrast with Theorem 4.5 from \cite{toeplitz-projection} which treats the case $p=\infty$,
the tuple $T$ is not assumed to be essentially normal here. But unfortunately, the methods of the present paper seem not to carry over to
the case $p=\infty$, at least not in a straightforward way. So it remains open 
whether the essential-normality condition can be dropped from the hypotheses in the case $p = \infty$.

Starting point for the present note was a recent paper of Hasegawa \cite{hasegawa} where $p$-essential commutants in the 
context of semicrossed products are studied. Among other results it was shown in \cite{hasegawa} that Theorem \ref{esscom_p}
holds in the particular case where $T = M_z \in B(H^2(\mathbb D))$ is the multiplication operator with the argument on the 
Hardy space of the unit disc and $A = A(\mathbb D)$ is the disc algebra. 
However, the methods in \cite{hasegawa} are limited by the availability of an underlying abelian group structure
needed to define suitable semicrossed products.

Under the condition that $W^*(U)$ is a maximal abelian von Neumann algebra, each abstract $T$-Toeplitz operator possesses
an $L^\infty(\mu)$-symbol \cite{toeplitz-projection}. In this case, the above theorem yields in fact a characterization of $\Tcal_a^{\text{ec},p}$ in terms
of Hankel operators. In analogy with the classical case, these are defined by the formula
\[
H_f = (1-P_\Hcal)\Psi_U(f)|\Hcal \in B(\Hcal, \Kcal \ominus \Hcal) \qquad (f \in L^\infty(\mu)).
\]

\begin{cor} \label{esscom_p_concrete}
Let $T\in B(\Hcal)^n$ be a regular $A$-isometry with minimal normal extension $U \in B(\Kcal)^n$ and scalar spectral measure
$\mu \in M^+_1(\partial_A)$. If $W^*(U)$ is a maximal abelian von Neumann algebra, then an operator
$S \in B(\Hcal)$ belongs to $\Tcal_a(T)^{\text{ec},p}$ if and only if $S = T_f + K$ with some $K \in \Scal_p(\Hcal)$ and 
a symbol $f \in L^\infty(\mu)$ having the property that $H_f \in \Scal_p(\Hcal, \Kcal \ominus \Hcal)$.
\end{cor}

If $(A|\partial_A,\partial_A,\mu)$ is regular in the sense of Aleksandrov, then $T = M_z \in B(H^2_A(\mu))^n$ is a
regular $A$-isometry with minimal normal extension $U = M_z \in B(L^2(\mu))^n$. Since $W^*(U) = \{M_f: f \in L^{\infty}(\mu)\}$
is maximal abelian, Corollary \ref{esscom_p_concrete} gives a complete description of the $p$-essential commutant
$\Tcal_a(T)^{ec,p}$ in this case. In particular, Corollary \ref{esscom_p_concrete} applies to Hardy spaces on strictly
pseudoconvex and bounded symmetric domains.

Similarly, Theorem \ref{xia} will follow as a particular case of a more general result for regular $A$-isometries.
\vspace{1cm}

\centerline{\textbf{\S2\, Analytic Toeplitz operators}}

Let us fix a regular $A$-isometry $T \in B(\Hcal)^n$ and let $U$, $\Psi_U$ and $\mu$ be defined as in the introduction.
Our proof of Theorem \ref{esscom_p} makes use of a Toeplitz projection $\Phi_T$ established in \cite{toeplitz-projection} by
Everard and the second author. We briefly recall some facts about $\Phi_T$ that will be used in the sequel. First of all,
$\Phi_T$ is a completely positive unital projection
\[
\Phi_T : B(\Hcal) \to B(\Hcal) \quad \text{with} \quad \Phi_T(B(\Hcal)) = \Tcal(T).
\]
The operators $\Phi_T(X)$ are obtained as weak$^*$-limits $\Phi_T(X) = \lim_{\alpha}\Phi_{T,k_{\alpha}}(X)$, where
$(\Phi_{T,k_{\alpha}})$ is a subnet of the sequence of mappings
$\Phi_{T, k_\alpha}: B(\Hcal) \to B(\Hcal)$ defined by
\[
\Phi_{T,k}(X) = \frac{1}{k^k} \sum_{1 \leq i_1, \ldots, i_k \leq k} 
T^*_{\theta_{k}^{i_k} \cdots \theta_{1}^{i_1}} X T_{ \theta_{1}^{i_1} \cdots \theta_{k}^{i_k}} \qquad (k \geq 1),
\]
and $(\theta_k)_{k\geq 1}$ is a fixed sequence of $\mu$-inner functions with $W^*((\theta_k)_{k\geq 1}) = L^\infty(\mu)$.

\textbf{Proof of Theorem \ref{esscom_p}.} Let $S \in \Tcal_a^{\text{ec},p}$ be given.
Since the Schatten-$p$-norm dominates the operator norm, a straightforward application of the closed graph theorem shows that the linear map
\[
C_S : H^\infty_A(\mu) \to \Scal_p(\Hcal), \quad g \mapsto [S, T_g]
\]
is norm-bounded. Using the definition of the Toeplitz projection $\Phi_T$ explained above,
one obtains that the operator $\Phi_T(S)-S$ is the weak$^*$-limit of a net consisting of convex combinations 
of operators of the form
\[
T_\theta^*[S,T_\theta] \qquad (\theta \in I_\mu).
\]
Since these convex combinations are bounded in the Schatten-$p$-norm by the norm of the operator $C_S$, the lower semicontinuity
of the Schatten-$p$-norm with respect to the weak operator topology \cite{hiai} (Proposition 2.11) implies that $\Phi_T(S)-S \in \Scal_p(\Hcal)$. But then
\[
S = \Phi_T(S) +(S- \Phi_T(S)) \in \Tcal(T) + \Scal_p(\Hcal),
\]
as was to be shown. \proofend

To prepare the proof of Corollary \ref{esscom_p_concrete} we formulate a series of elementary lemmas related to Hankel operators.
To begin with observe that, for every $\mu$-inner function $\theta \in I_\mu$, the operator $\Psi_U(\theta) \in B(\Kcal)$ is unitary 
and leaves $\Hcal$ invariant.

\begin{lem} \label{hank}
For each $\mu$-inner function $\theta \in I_\mu$, the Hankel operator $H_{\overline{\theta}} \in B(\Hcal, \Kcal\ominus \Hcal)$ is a partial isometry
with
\[
\ker H_{\overline{\theta}} = \Psi_U(\theta)\Hcal \qquad \text{and} \qquad \ran H_{\overline{\theta}} = (\Kcal \ominus \Hcal) \ominus
\big(\Psi_U(\overline{\theta}) (\Kcal \ominus \Hcal)\big).
\]
\end{lem}
\proof
Obviously, $H_{\overline{\theta}} \Psi_U(\theta) \Hcal = P_{\Kcal \ominus \Hcal} \Psi_U(\overline{\theta})\Psi_U(\theta) \Hcal = \{0\}$. Since
on the other hand,
$\Psi_U(\overline{\theta}) (\Kcal \ominus \Psi_U(\theta)\Hcal) = \Kcal \ominus \Hcal$, it follows that $H_{\overline{\theta}} = \Psi_U(\overline{\theta})$
on $\Hcal \ominus \Psi_U(\theta)\Hcal$. Hence $H_{\overline{\theta}}$ is a partial isometry with $\ker H_{\overline{\theta}} = \Psi_U(\theta)\Hcal$.
The observation that 
\[
H_{\overline{\theta}} ( \Hcal \ominus \Psi_U(\theta)\Hcal) = (\Psi_U(\overline{\theta})\Hcal) \ominus \Hcal
= (\Kcal \ominus \Hcal) \ominus \big(\Psi_U(\overline{\theta}) (\Kcal \ominus \Hcal) \big)
\]
completes the proof. \proofend

By Lemma \ref{hank}, for each $\mu$-inner function $\theta \in I_\mu$, the operator
\[
P_\theta = H_{\overline{\theta}} H_{\overline{\theta}}^* \in B(\Kcal \ominus \Hcal)
\]
is the orthogonal projection of $\Kcal \ominus \Hcal$ onto the space 
\[
\Kcal_\theta = (\Kcal \ominus \Hcal) \ominus \big(\Psi_U(\overline{\theta})(\Kcal\ominus \Hcal)\big).
\]
\begin{lem} \label{span}
With the notations from above, we have $\Kcal \ominus \Hcal = \bigvee \left( \Kcal_\theta : \theta \in I_\mu \right)$.
\end{lem}
\proof The orthogonal complement of the space on the right with respect to $\Kcal \ominus \Hcal$ is given by
\[
\bigcap_{\theta \in I_\mu} \Psi_U(\overline{\theta}) (\Kcal \ominus \Hcal)
= \bigcap_{\theta \in I_\mu} \Kcal \ominus (\Psi_U(\overline{\theta})\Hcal)
= \Kcal \ominus \bigvee_{\theta \in I_\mu} \Psi_U(\overline{\theta})\Hcal.
\]
The closed linear span $\bigvee_{\theta \in I_\mu} \Psi_U(\overline{\theta})\Hcal$ contains $\Hcal$ and is invariant under the von Neumann
algebra $W^*(\Psi_U(I_\mu)) = \Psi_U(W^*(I_\mu)) = W^*(U)$. Thus the assertion follows from the minimality of $U$ as a normal extension
of $T$. \proofend

There is a canonical partial order on the set $I_\mu$ of all $\mu$-inner functions. By definition $\theta_1 \leq \theta_2$ if there is a $\mu$-inner
function $\theta \in I_\mu$ with $\theta_2 = \theta_1 \theta$. Obviously, the partially ordered set $(I_\mu, \leq)$ is directed upwards.
If $\theta_1 \leq \theta_2$ and $\theta_2 = \theta_1 \theta$ as above, then
\[
\Psi_U(\overline{\theta}_2) (\Kcal \ominus \Hcal) = \Psi_U(\overline{\theta}_1) \big(\Psi_U(\overline{\theta})(\Kcal \ominus \Hcal)\big)
\subset \Psi_U(\overline{\theta}_1) (\Kcal \ominus \Hcal)
\]
and hence $P_{\theta_1} \leq P_{\theta_2}$.

\begin{lem}\label{SOT-identity}
The net $(P_\theta)_{\theta \in I_\mu}$ converges strongly to the identity operator $I_{\Kcal\ominus \Hcal}$ on $\Kcal\ominus \Hcal$.
\end{lem}
\proof By Lemma \ref{span}, the space $\Kcal \ominus \Hcal$ is the closed linear span of the set
\[
M= \{ P_\theta h : \theta \in I_\mu \text{ and } h \in \Kcal \ominus \Hcal \}.
\]
Since the net $(P_\theta)_{\theta \in I_\mu}$ converges
pointwise to the identity operator on $M$, it converges pointwise to the identity operator on the closed linear span $\Kcal \ominus \Hcal$ of $M$. \proofend

With these preparations we can prove the claimed result in the maximal abelian case.

\textbf{Proof of Corollary \ref{esscom_p_concrete}.}
Suppose first that $S = T_f +K$ with $K \in \Scal_p(\Hcal)$ and a symbol $f \in L^\infty(\mu)$ such that $H_f \in \Scal_p(\Hcal, \Kcal\ominus \Hcal)$.
To prove that $S \in \Tcal_a(T)^{\text{ec},p}$, it suffices to observe that, for $g \in H^\infty_A(\mu)$ and $f \in L^\infty(\mu)$, we have
the identity 
\[
[T_f, T_g] = T_{fg} - T_gT_f = H_{\overline{g}}^* H_f.
\]
Conversely, suppose that $S\in \Tcal_a(T)^{\text{ec},p}$. 
By Theorem \ref{esscom_p} we know that $S = X+K$ with $X \in \Tcal(T)$ and $K \in \Scal_p(\Hcal)$. Since $W^*(U)$ is maximal abelian
it follows that $\Tcal(T) = P_\Hcal W^*(U)'|\Hcal = P_\Hcal W^*(U)|\Hcal$. Hence there is a function $f \in L^\infty(\mu)$
with $X = T_f$. Using the formula for $[T_f, T_g]$ from the first part of the proof, we find that
$H_{\overline{\theta}}^* H_f \in \Scal_p(\Hcal)$ for all $\mu$-inner functions $\theta \in I_\mu$. By Lemma \ref{SOT-identity}, the Hankel operator
\[
H_f = \text{SOT}-\lim_{\theta \in I_\mu} P_\theta H_f = \text{SOT}-\lim_{\theta \in I_\mu} H_{\overline{\theta}} (H_{\overline{\theta}}^*H_f)
\]
is the strong limit of a net in $\Scal_p(\Hcal, \Kcal \ominus \Hcal)$. With the notations from the proof of Theorem \ref{esscom_p}
we obtain that
\[
\| P_\theta H_f \|_p \leq \| H_{\overline{\theta}} \| \, \| H_{\overline{\theta}}^* H_f \|_p =
\| H_{\overline{\theta}} \| \, \| [T_f, T_\theta]\|_p \leq \| C_{T_f} \|
\]
for all $\theta \in I_\mu$. Using again the semicontinuity of the Schatten-$p$-norm with respect to the
weak operator topology, we conclude that $H_f \in \Scal_p(\Hcal, \Kcal \ominus \Hcal)$.
\proofend

\textbf{Proof of Theorem \ref{esscom_p_ball_polydisc}.}
With $D = \B_n$ (or $D = \D^n$, respectively) we apply Corollary \ref{esscom_p_concrete} to the canonical
probability measure $\sigma$
on $\partial \B_n$ ($\T^n$, respectively) and $T = M_z \in B(H^2_{A(D)})^n$.
In both cases, $H^2_{A(D)}$ is the classical Hardy space on the corresponding Shilov boundary of $D$.
If $S \in \Tcal_a^{\text{ec},p}$, then the cited corollary implies that $S = T_f + K$ with $K \in \Scal_p$ and
\[
H_f \in \Scal_{p}(H^2_{A(D)}, L^2(\sigma) \ominus H^2_{A(D)}).
\]
By our assumption that $1\leq p \leq 2n$ in the ball-case, Theorem 1.5 from Xia and Fang \cite{fang-xia} says that
$H_f \in \Scal_{2n}$ is the zero operator. 
In the polydisc-case a theorem of Cotlar and Sadoski \cite{cotlar} asserts that there are no non-zero compact Hankel
operators.
In both cases, this implies that the symbol $f \in L^\infty(\sigma)$ actually belongs to $f \in H^\infty(\sigma)$.
Since the other direction is trivial, the proof is complete.
\proofend

Let $T \in B(\Hcal)^n$ be a regular $A$-isometry. The methods used above can be used in a similar way to describe the
commutant of $\Tcal_a(T)$ modulo the set $\Scal_0(\Hcal)$ of finite-rank operators on $\Hcal$, at least in the case where
$T$ has empty point spectrum. By definition the point spectrum of $T$ is
\[
\sigma_p(T) = \{\lambda \in \C^n: \bigcap^n_{i=1} \ker(\lambda_i - T_i) \neq \{0\} \}.
\]
If $\sigma_p(T) = \varnothing$, then the scalar spectral measure $\mu \in M^+_1(\partial_A)$ of the minimal normal extension 
$U \in B(\Kcal)^n$ is continuous, and Aleksandrov's results on the abstract inner function problem imply that there is a
sequence $(J_k)_{k\geq 1}$ of isometries in $\Tcal_a(T)$ with $w^*$-$\lim_{k \rightarrow \infty} J_k = 0$ in $B(\Hcal)$.
Exactly as in the particular case of spherical isometries (Theorem 3.5 in \cite{inner-spherical}), the sequence $(J_k)_{k\geq 1}$
can be used to prove the first part of the following result.

\begin{thm}\label{finite}
Let $T \in B(\Hcal)^n$ be a regular $A$-isometry with minimal normal extension $U \in B(\Kcal)^n$ and scalar spectral measure
$\mu \in M^+_1(\partial_A)$. Suppose that $\sigma_p(T) = \varnothing$. Then every operator $S \in \Tcal_a(T)^{ec,0}$ has the form
\[
S = X + K \qquad \mbox{with} \; X \in \Tcal(T) \; \mbox{and} \;  K \in \Scal_0(\Hcal).
\]
If in additdion $W^*(U)$ is maximal abelian, then
\[
\Tcal_a(T)^{ec,0} = \{T_f + K: \;  K \in \Scal_0(\Hcal) \; \mbox{and} \; f \in L^{\infty}(\mu) \; \mbox{with} \; H_f \in \Scal_0(\Hcal,\Kcal \ominus \Hcal)\}.
\]
\end{thm}

\proof
For the proof of the first part, we refer to Theorem 3.5 in \cite{inner-spherical}. Suppose that $W^*(U)$ is maximal abelian.
Exactly as in the proof of Corollary \ref{esscom_p_concrete} the formula
\[
[T_f,T_g] = H^*_{\overline{g}} H_f \qquad (g\in H_A^{\infty}(\mu), f \in L^{\infty}(\mu))
\]
shows that $T_f + K \in \Tcal_a(T)^{ec,0}$ for $K \in \Scal_0(\Hcal)$ and every function $f \in L^{\infty}(\mu)$ with the property that $H_f$ is a finite-rank
operator. Conversely, let $f \in L^{\infty}(\mu)$ be a symbol such that $[T_f,T_g]$ is a finite-rank operator for all $g \in H_A^{\infty}(\mu)$. 
Since $\Tcal_a(T) \subset B(\Hcal)$ is closed (Corollary 3.6 in \cite{toeplitz-projection}), it follows from part (b) of Lemma 3.4 in
\cite{inner-spherical} that there is a constant $M > 0$ such that ${\rm rank}([T_f,T_g]) \leq M$ for all $g \in H_A^{\infty}(\mu)$. As in the proof of
Corollary \ref{esscom_p_concrete} we find that
\[
H_f = {\rm SOT}-\lim_{\theta \in I_{\mu}} H_{\overline{\theta}}(H^*_{\overline{\theta}} H_f) = {\rm SOT}-\lim_{\theta \in I_{\mu}} H_{\overline{\theta}}[T_f,T_{\theta}].
\]
Hence part (a) of Lemma 3.4 from \cite{inner-spherical} shows that ${\rm rank} \, H_f \leq M$. Since by the first part of Theorem \ref{finite} and the
hypothesis that $W^*(U)$ is maximal abelian every operator $S \in \Tcal_a(T)^{ec,0}$ has the form $S = T_f + K$ with $f \in L^{\infty}(\mu)$ and
$K \in \Scal_0(\Hcal)$, the above observations complete the proof.
\proofend

As particular cases, Theorem \ref{finite} contains the results of Guo and Wang \cite{guo-wang} on operators commuting  with analytic 
Toeplitz operators modulo funite-rank operators on the unit ball and the unit polydisc. To obtain the concrete form of the results from
\cite{guo-wang} it suffices to use Kronecker's theorem on finite-rank Hankel operators on $\mathbb D$, and the results on Schatten class
Hankel operators on $\B_n$ and $\mathbb D^n$ cited in the proof of Theorem \ref{esscom_p_ball_polydisc}.

\vspace{1cm}

\centerline{\textbf{\S3 \, On Schatten class perturbations of Toeplitz operators}}

Let $D\subset \mathbb C^n$ be a bounded domain. In this section we specialize to the case that $T \in B(\Hcal)^n$ is a regular $A$-isometry with respect to the domain algebra
\[
A=A(D)=\lbrace f \in C(\overline{D}):\ f|D \in \mathcal O(D)\rbrace.
\]
Furthermore, we make the additional assumption that the scalar spectral measure $\mu \in M^+_1(\partial_A)$ associated with $T$ is a faithful Henkin measure.  
The condition that 
$\mu \in M^+_1(\partial_A)$ is a Henkin measure means by definition that there is a contractive and weak$ ^\ast$ continuous algebra homomorphism
\[
r_\mu:H^\infty(D)\rightarrow L^\infty(\mu)
\]
which extends the canonical map $A(D)\rightarrow L^\infty(\mu),\ f \mapsto [f|_{\partial_A}]$. Here the weak$ ^\ast$ closed subspace $H^\infty(D)\subset L^\infty(D)$ is regarded as the dual space of
the quotient space $L^1(D)/^\perp H^\infty(D)$. If $r_\mu$ is even isometric, then $\mu$ is said to be a faithful Henkin measure. In this case the range 
$\ran r_\mu\subset L^\infty(\mu)$ is 
weak$ ^\ast$ closed and
\[
r_\mu:H^\infty(D)\rightarrow \ran r_\mu
\]
is a dual algebra isomorphism, that is, an isometric isomorphism of Banach algebras and a weak$ ^\ast$ homeomorphism. In particular, the range of $r_\mu$ contains the weak$ ^\ast$ closed subalgebra
\[
H^\infty_A(\mu)=\overline{A|_{\partial_A}}^{w^\ast}\subset L^\infty(\mu).
\]
To simplify the notation we often write $f^\ast=r_\mu(f)\in L^\infty(\mu)$ for the boundary value of a function $f \in H^\infty(D)$.
As before, we denote by
\[
I_\mu=\lbrace \theta \in H^\infty_A(\mu):\ |\theta|=1\ \mu\mbox{-almost everywhere on }\partial_A\rbrace
\]
the set of all $\mu$-inner functions. We call the elements of
\[
I_D=\lbrace \theta \in H^\infty(D):\ r_\mu(\theta)\in I_\mu\rbrace
\]
the inner functions on $D$ (with respect to $\mu$) and write
\[
B_{H^\infty(D)}=\lbrace f \in H^\infty(D):\ \| f \|_D \leq 1\rbrace
\]
for the closed unit ball of $H^\infty(D)$. A well known and standard argument (cf. Lemma 14.1.6 in \cite{dales}) shows that a sequence $(f_k)_{k \in \mathbb N}$ in $H^\infty(D)$ is weak$ ^\ast$ convergent to a function
$f \in H^\infty(D)$ if and only if $(\| f_k\|_D)_{k \in \mathbb N}$ is bounded and $(f_k)_{k \in \mathbb N}$ converges to $f$ pointwise, or equivalently, uniformly on all compact
subsets of $D$.

The condition that $(A|\partial_A,\partial_A,\mu)$ is regular implies that the space $L^{\infty}(\mu)$ is generated as a von Neumann algebra
by $I_{\mu}$ (Corollary 2.5 in \cite{inner-spherical}). Since $I_{\mu} = r_{\mu}(I_D)$ and since $r_{\mu}: H^{\infty}(D) \rightarrow L^{\infty}(\mu)$ is
injective, it follows that $I_D$ contains at least one non-constant inner function $\theta$. But then $(\theta^k)_{k \in \mathbb N}$ is a weak$^*$ zero 
sequence in $H^{\infty}(D)$ and $(r_{\mu}(\theta^k))_{k \in \mathbb N}$ is a weak$^*$ zero sequence in $L^{\infty}(\mu)$. As an elementary consequence we
obtain that the measure $\mu$ is necessarily continuous. Indeed, if there were a point $z \in \partial_A$ with $\mu(\{z\}) > 0$, then we could
conclude that $|\theta^*(z)| = 1$, which is impossible since
\[
\theta^*(z)^k \mu(\{z\}) = \langle \chi_{\{z\}}, r_{\mu}(\theta^k)\rangle_{(L^1(\mu),L^{\infty}(\mu))} \stackrel{k}{\longrightarrow} 0.
\]

We need conditions which ensure that the infinite product a given sequence $(\theta_j)_{j \in \mathbb N}$ of inner functions on $D$ is an inner function again.

\begin{lem}\label{prod1}
Let $w \in D$ be arbitrary and let $(u_j)_{ j \in \mathbb N}$ be a sequence in $B_{H^\infty(D)}$ such that $u_j(w)\neq 0$ for all $j$ and such that
\[
\sum^\infty_{j=0}|1-u_j(w)|< \infty.
\]
Then the sequence $(v_N)_{N \in \mathbb N}=\big(\prod^N_{j=0}u_j\big)_{N \in \mathbb N}$ is weak$^\ast$ convergent to a function $v \in B_{H^\infty(D)}$ with $v(w)\neq 0$.
\end{lem}

\proof
The hypothesis of the lemma implies that the limit
\[
\lim_{N\rightarrow \infty}v_N(w)\in \mathbb C \setminus \lbrace 0 \rbrace
\]
exists. By Montel's theorem there is a subsequence $(v_{N_j})_{j\in \mathbb N}$ of $(v_N)_{N\in \mathbb N}$ which is weak$ ^\ast$ convergent to some function $v \in B_{H^\infty(D)}$.
Let $(v_{\tilde{N}_j})_{j \in \mathbb N}$ be another weak$^\ast$ convergent subsequence of $(v_N)_{N\in \mathbb N}$ with weak$ ^\ast$ limit $\tilde{v}$. Then
\[
| v(z)|=\lim_{N\rightarrow \infty}|v_N(z)|=|\tilde{v}(z)|\quad (z \in D)
\]
and hence $|v/\tilde{v}|=1$ on the complement of the zero set $Z(\tilde{v})$ of  $\tilde{v}$ in $D$. By Riemann's extension theorem and the open mapping principle from several complex variables
it follows that $v/\tilde{v}\equiv e^{i\theta}$ on $D\setminus Z(\tilde{v})$ for some real number $\theta \in \mathbb R$. Because of $v(w)=\tilde{v}(w)\neq 0$ we obtain that
$v=\tilde{v}$. The above arguments show that each subsequence of $(v_N)_{N\in \mathbb N}$ has a subsequence which is weak$^\ast$ convergent to $v$. But then it easily follows that the
sequence $(v_N)_{N\in \mathbb N}$ itself is weak$^\ast$ convergent to $v$.
\proofend

For $g \in L^\infty(\mu)$, we denote by
\[
M_g:L^2(\mu)\rightarrow L^2(\mu),\ h \mapsto gh
\]
the induced multiplication operator.

\begin{thm}\label{1 to 6}
For a given sequence $(\theta_k)_{k \in \mathbb N}$ in $B_{H^\infty(D)}$, the following conditions are equivalent:
	\begin{enumerate}[(i)]
		\item there is a point $w \in D$ such that $\lim_{k \rightarrow \infty}\theta_k(w)=1$,
		\item $\lim_{k\rightarrow \infty}\ \|\theta^\ast_k - 1\|_{L^2(\mu)}=0$,
		\item ${\rm SOT}-\lim_{k \rightarrow \infty} \ M_{\theta^\ast_k}=1_{L^2}(\mu)$,
		\item ${\rm WOT}-\lim_{k \rightarrow \infty}\ M_{\theta^\ast_k}=1_{L^2(\mu)}$,
		\item $w^\ast$-$\lim_{k \rightarrow \infty}\theta^\ast_k=1$ with respect to the weak$ ^\ast$ topology of $L^\infty(\mu)=L^1(\mu)^\prime,$
		\item $\lim_{k \rightarrow \infty}\int_{\partial_A}\theta^\ast_k\ d\mu=1$.
	\end{enumerate}
\end{thm}

\proof
Suppose that there is a point $w \in D$ with $\lim_{k \rightarrow \infty}\theta_k(w)=1$. Then the limit $\theta$ of any weak$^\ast$ convergent subsequence of $(\theta_k)_{k \in \mathbb N}$ is a function
$\theta \in B_{H^\infty(D)}$ with $\theta(w)=1$. By the maximum modulus principle it follows that necessarily $\theta \equiv 1$. As in the proof of Lemma \ref{prod1} we may conclude that $(\theta_k)_{k \in \mathbb N}$ is 
weak$^\ast$ convergent to the constant function $1$. The weak$^\ast$ continuity of the boundary value ${\rm map}\ r_\mu$ implies that $w^\ast-\lim_{k \rightarrow \infty}\theta^\ast_k=1$ in
$L^\infty(\mu)=L^1(\mu)^\prime$. In particular,
\[
\int\limits_{\partial_A}\theta^\ast_k\ d\mu=\langle 1,\theta^\ast_k\rangle \stackrel{k}{\longrightarrow} \int\limits_{\partial_A}1\ d\mu=1.
\]
But then we obtain that
\[
\int\limits_{\partial_A}|1-\theta^\ast_k|^2d\mu\leq \int\limits_{\partial_A}2(1-{\rm Re}\ \theta^\ast_k)d\mu\stackrel{k}{\longrightarrow} 0.
\]
Let us assume that $\lim_{k\rightarrow \infty}\theta^\ast_k=1$ in $L^2(\mu)$. Then there is a subsequence $(\theta^\ast_{k_j})_{j \in \mathbb N}$ of $(\theta^\ast_k)_{k \in \mathbb N}$ which converges $\mu$-almost
everywhere on $\partial_A$ to the constant function $1$. By the dominated convergence theorem
\[
{\rm SOT}-\lim_{j \rightarrow \infty}M_{\theta^\ast_{k_j}}=1_{L^2(\mu)}.
\]
By applying this argument to each subsequence of $(\theta^\ast_k)_{k \in \mathbb N}$ we find that
\[
{\rm SOT}-\lim_{k \rightarrow \infty}M_{\theta^\ast_k}=1_{L^2(\mu)}.
\]
To conclude the proof, let us suppose that ${\rm WOT}-\lim_{k \rightarrow \infty}M_{\theta^\ast_k}=1_{L^2(\mu)}$. Since each function in $L^1(\mu)$ is the product of two functions in $L^2(\mu)$, it follows that
\[
w^\ast-\lim_{k \rightarrow \infty}\theta^\ast_k=1 \mbox{ in } L^\infty(\mu)=L^1(\mu)^\prime.
\]
But then the fact that $r_\mu:H^\infty(D)\rightarrow \ran r_\mu$ is a weak$ ^\ast$ homeomorphism implies that $w^\ast-\lim_{k\rightarrow \infty}\theta_k=1$ in $H^\infty(D)$ and hence that $(\theta_k)\stackrel{k}{\longrightarrow} 1$ uniformly on all compact subsets of $D$.
\proofend

An inspection of the above proof shows that the equivalence of the conditions (ii)-(vi) holds more generally with $(\theta^\ast_k)_{k \in \mathbb N}$ replaced by an arbitrary sequence $(f_k)_{k \in \mathbb N}$ in
$B_{L^\infty(\mu)}$.

\begin{cor}\label{prod2}
Let $w \in D$ and let $(u_j)_{j \in \mathbb N}$ be a sequence in $B_{H^\infty(D)}$ with $u_j(w)\neq 0$ for all $j \in \mathbb N$ and
\[
\sum^\infty_{j=0}| 1-u_j(w)|< \infty.
\]
Then the sequence $(v_N)_{N \in \mathbb N}=\big(\prod^N_{j=0}u_j\big)_{N \in \mathbb N}$ of partial products is weak$ ^\ast$ convergent in $H^\infty(D)$ to some function $v \in H^\infty(D)$. Furthermore,
\[
\lim_{N\rightarrow \infty}\| v^\ast_N-v^\ast\|_{L^2(\mu)}=0
\]
and $M_{v^\ast}={\rm SOT}-\lim_{N \rightarrow \infty}M_{v^\ast_N}$ in $B(L^2(\mu))$. If each  $u_j$ is an inner function on $D$, then so is $v$.
\end{cor}

\proof
By Lemma \ref{prod1} the limit $v=w^\ast-\lim_{N\rightarrow \infty}v^\ast_N$ exists in $H^\infty(D)$. By applying Lemma \ref{prod1} to the sequences $(u_j)_{j\geq N+1}$ we obtain in the same way that the products
\[
\theta_N=\prod^\infty_{j=N+1}u_j\in B_{H^\infty(D)}
\]
converge uniformly on all compact subsets of $D$, or equivalently, with respect to the weak$ ^\ast$ topology of $H^\infty(D)$. The estimates
\[
|\theta_N(w)-1|=\lim_{k\rightarrow \infty}\left|\left(\prod^k_{j=N+1}1+(u_j(w)-1)\right)-1\right|
\]
\[
\leq \lim_{k \rightarrow \infty}\left[\left(\prod^k_{j=N+1}1+|u_j(w)-1|\right)-1\right]\leq \lim_{k\rightarrow \infty}\left(e^{\sum^k_{j=N+1}|u_j(w)-1|}-1\right)
\]
imply that $\lim_{N \rightarrow \infty}\theta_N(w)=1$. But then Theorem \ref{1 to 6} yields that
\[
v^\ast-v^\ast_N=v^\ast_N(\theta^\ast_N-1)\stackrel{N}{\longrightarrow}0 \mbox{ in }L^2(\mu)
\]
and that
\[
M_{v^\ast}-M_{v^\ast_N}=M_{v^\ast_N}(M_{\theta^\ast_N}-1_{L^2(\mu)})\stackrel{N}{\longrightarrow} 0
\]
in $B(L^2(\mu))$ with respect to the strong operator topology. Since $(v^\ast_N)_{N\in \mathbb N}$ possesses a subsequence which converges pointwise $\mu$-almost everywhere to $v^\ast$, also the last assertion follows.
\proofend

Let $(S_j)_{j\in \mathbb N}$ be a sequence of bounded operators $S_j\in B(\Hcal)$. If ${\rm SOT}-\lim_{j\rightarrow \infty}S_j=0$, then a standard compactness argument allows one to conclude that
\[
\lim_{j\rightarrow \infty}\| S_j K\|=0 \mbox{   and   }\lim_{j\rightarrow \infty}\| KS^\ast_j\|=\lim_{j \rightarrow \infty}\| S_j K^\ast\| =0
\]
for each compact operator $K \in K(\Hcal)$ on $\Hcal$. Using the well known fact that each Schatten class operator $S\in \mathcal S_p(\Hcal)\quad (1\leq p< \infty)$ can be written as a product $S=K\tilde{S}$ of
a compact operator $K\in K(\Hcal)$ and an operator $\tilde{S}\in \mathcal S_p(\Hcal)$, one obtains the following elementary consequence.

\begin{lem}\label{SOT}
Let $(S_j)_{j \in \mathbb N}$ be a sequence in $B(\Hcal)$ such that
\[
{\rm SOT}-\lim_{j\rightarrow \infty}S_j=0.
\]
Then, for $1\leq p \leq \infty$ and $S \in \mathcal S_p(\Hcal)$, we have
\[
\lim_{j \rightarrow \infty}\| S_jS\|_p=\lim_{j \rightarrow \infty}\| SS^\ast_j\|_p=0.
\]
\end{lem}

\proof
For $p<\infty$, write $S=K\tilde{S}$ with $K\in K(\Hcal)$ and $\tilde{S}\in \mathcal S_p(\Hcal)$. Then
\[
\| S_jS\|\leq \| S_jK\|\ \|\tilde{S}\|_p\stackrel{j}{\longrightarrow}0
\]
and $\| SS^\ast_j\|_p=\|S_jS^\ast\|_p\stackrel{j}{\longrightarrow} 0$, since also $S^\ast \in \mathcal S_p(\Hcal)$.
\proofend

We need some additional continuity properties of the Toeplitz calculus.

\begin{lem}\label{biSOT}
Let $(f_k)_{k \in \mathbb N}$ be a bounded sequence in $L^\infty(\mu)$ such that
\[
\lim_{k \rightarrow \infty}\| f_k - f\|_{L^2(\mu)}=0
\]
for some function $f \in L^\infty(\mu)$. Then
\[
{\rm SOT}-\lim_{k \rightarrow \infty}T_{f_k}=T_f \quad \mbox{ and }\quad {\rm SOT}-\lim_{k\rightarrow \infty}T^\ast_{f_k}=T^\ast_f.
\]
\end{lem}

\proof
Each subsequence of $(f_k)_{k\in \mathbb N}$ possesses a subsequence $(g_k)_{k\in \mathbb N}$ such that $(g_k)\stackrel{k}{\longrightarrow} f$ pointwise $\mu$-almost everywhere
on $\partial_A$. By Lemma 3.4 in \cite{esscom} the sequence $(T_{g_k})_{k \in \mathbb N}$ is {\rm SOT}-convergent to $T_f$. But then an elementary argument shows that ${\rm SOT}-\lim_{k\rightarrow \infty}
T_{f_k}=T_f$. Since $(\overline{f}_k)_{k\in \mathbb N}$ and $\overline{f}$ satisfy the same hypotheses as $(f_k)_{k\in \mathbb N}$ and $f$, also the second assertion follows.
\proofend

Up to now we only used  Toeplitz operators with symbol $f \in L^\infty(\mu)$. To simplify the notation, we define $T_f=T_{r_\mu(f)}$ for $f \in H^\infty(D)$.

\begin{cor}\label{xia1to3}
Let $(\theta_k)_{k\in \mathbb N}$ be a sequence in $I_D$ such that there is a point $w \in D$ with $\lim_{k\rightarrow \infty}(1-|\theta_k(w)|)=0$. Then, for $1 \leq p\leq \infty$, we have
	\begin{enumerate}[(i)]
		\item ${\rm SOT}-\lim_{k\rightarrow \infty}T^\ast_{\theta_k}X T_{\theta_k}=X$ for all $X \in B(\Hcal)$,
		\item $\lim_{k \rightarrow \infty}\| T^\ast_{\theta_k}X T_{\theta_k}-X\|_p=0$ for all $X\in \mathcal S_p(\Hcal)$,
		\item if $X\in B(\Hcal)$ and $u \in H^\infty(D)$ with $T^\ast_uXT_u-X\in \mathcal S_p(\Hcal)$, then
		\[
		\lim_{k \rightarrow \infty}\| T^\ast_u(T^\ast_{\theta_k}XT_{\theta_k}-X)T_u-(T^\ast_{\theta_k}XT_{\theta_k}-X)\|_p=0.
		\]
	\end{enumerate}
\end{cor}

\proof
Since $T^\ast_{\theta_k}XT_{\theta_k}=T^\ast_{\epsilon\theta_k}XT_{\epsilon\theta_k}$ for each complex number $\epsilon$ with $|\epsilon|=1$, we may assume that $\theta_k(w)\geq 0$ for all $k$.
Then part (i) follows immediately from Theorem \ref{1 to 6} and Lemma \ref{biSOT}. For $X \in \mathcal S_p(\Hcal)$, a combination of Lemma \ref{SOT} and Lemma 
\ref{biSOT} yields that
\begin{eqnarray*}
\| T^\ast_{\theta_k}XT_{\theta_k}-X\|_p&=& \| T^\ast_{\theta_k}(XT_{\theta_k}- T_{\theta_k}X)\|_p\\
&\leq& \| X(T_{\theta_k}-1)\|_p+\|(1-T_{\theta_k})X\|_p \stackrel{k}{\longrightarrow} 0
\end{eqnarray*}
for $1 \leq p\leq \infty$. Fix an operator $X \in B(\Hcal)$ and a function $u \in H^\infty(D)$ such that $T^\ast_u XT_u-X\in \mathcal S_p(\Hcal)$ for some $p \in [1,\infty]$. Then
\begin{eqnarray*}
\| T^\ast_u(T^\ast_{\theta_k}XT_{\theta_k}-X)T_u&-&(T^\ast_{\theta_k}XT_{\theta_k}-X)\|_p\\
=\| T^\ast_{\theta_k}(T^\ast_uXT_u-X)T_{\theta_k}&-&(T^\ast_uXT_u-X)\|_p \stackrel{k}{\longrightarrow}0
\end{eqnarray*}
by part (ii).
\proofend

Using the weak$^*$-continuity of the $L^{\infty}(\mu)$-functional calculus of the minimal normal extension $U$ of $T$ 
we obtain the following variant of part (i) of Corollary \ref{xia1to3}.

\begin{lem}\label{compact}
Let $(v_k)_{k \in \mathbb N}$ be a sequence in $L^{\infty}(\mu)$ with
$w^*-\lim_{k\rightarrow \infty} = v \in L^{\infty}(\mu)$. Then, for each compact operator $K \in K(\Hcal)$,
\[
{\rm WOT}-\lim_{k\rightarrow \infty} T^*_{v_k}KT_{v_k} = T^*_v K T_v.
\]
\end{lem}

\proof Since the functional calculus map $\Psi_U: L^{\infty}(\mu) \rightarrow B(\Kcal)$ of the minimal normal
extension $U \in B(\Kcal)^n$ of $T$ and the compression mapping $B(\Kcal) \rightarrow B(\Hcal), X \mapsto P_{\Hcal}X|\Hcal,$
are weak$^*$ continuous, also the mapping
\[
L^{\infty}(\mu) \rightarrow B(\Hcal), f \mapsto T_f = P_{\Hcal} \Psi_U(f)|\Hcal
\]
is weak$^*$ continuous. Hence, under the hypothesis of the lemma, we find that
\[
{\rm SOT}-\lim_{k \rightarrow \infty} K T_{v_k} = K T_v
\]
for each compact operator $K \in K(\Hcal)$. But then
\[
T^*_{v_k}KT_{v_k} - T^*_v K T_v = T^*_{v_k}K(T_{v_k} - T_v) + (T^*_{v_k} - T^*_v)KT_v \stackrel{k}{\rightarrow} 0
\]
converges to zero in the weak operator topology.
\proofend

As mentioned before, the abstract $T$-Toeplitz operators are precisely the operators of the form
\[
T_Y = P_{\Hcal}Y|\Hcal \in B(\Hcal)
\]
with $Y \in (U)'$. If $S = T_Y + K$ with $Y \in (U)'$ and $K \in \mathcal S_p(\Hcal)$, then using the fact that
$T^*_{\theta} T_Y T_{\theta} = T_Y$ for all $\theta \in I_{\mu}$, we find that
\[
T^*_{\theta} S T_{\theta} - S \in \mathcal S_p(\Hcal)
\]
for $\theta \in I_{\mu}$. Our aim is to show that conversely each operator $S\in B(\Hcal)$ satisfying this
condition is a Schatten-$p$-class perturbation of a $T$-Toeplitz operator. The proof is based on an
observation which in the particular case of Toeplitz operators on the Hardy space of the unit disc $H^2(\mathbb D)$ is due to
Xia \cite{xia-characterization}.

Let $1 \leq p \leq \infty$. Suppose that $S = T_Y + K$ with $Y \in (U)'$ and $K \in \mathcal S_p(\Hcal)$. Then, for
each weak$^*$ convergent sequence $(\theta_k)_{k \in \mathbb N}$ in $I_{\mu}$ such that its weak$^*$ limit
$w^*-\lim_{k \rightarrow \infty} \theta_k = \alpha \in \mathbb C$ is a constant function, the limit
\[
X = {\rm WOT}-\lim_{k \rightarrow \infty} T^*_{\theta_k} S T_{\theta_k} = 
{\rm WOT}-\lim_{k \rightarrow \infty} T^*_{\theta_k} (T_Y + K) T_{\theta_k} = T_Y + |\alpha|^2 K
\]
exists by Lemma \ref{compact} and satisfies
\[
K = S - T_Y = S -X + |\alpha|^2 K = \frac{S - X}{1 - |\alpha|^2}
\]
as well as
\[
T_Y = S - K = \big( 1 - \frac{1}{1 - |\alpha|^2}\big) S + \frac{X}{1 - |\alpha|^2} = \frac{1}{1 - |\alpha|^2}\big(X - |\alpha|^2 S\big),
\]
provided that $|\alpha| \neq 1$.

Conversely, suppose that $S \in B(\Hcal)$ is an operator such that
\[
T^*_{\theta} S T_{\theta} - S \in \mathcal S_p(\Hcal)
\]
for all $\theta \in I_{\mu}$. Since $T$ is a regular $A$-isometry with a continuous scalar spectral measure $\mu$, Aleksandrov's
results (Corollary 29 in \cite{aleks}) on the existence of $\mu$-inner functions allow us to choose, for any given complex number $\alpha \in \mathbb C$ with
$|\alpha| < 1$, a sequence $(\theta_k)_{k \in \mathbb N}$ in $I_{\mu}$ such that
\[
w^*-\lim_{k \rightarrow \infty} \theta_k = \alpha \; {\rm in} \; L^{\infty}(\mu).
\]
By passing to a subsequence one can achieve that in addition the limit
\[
X = {\rm WOT}-\lim_{k \rightarrow \infty} T^*_{\theta_k} S T_{\theta_k} \in B(\Hcal)
\]
exists. Using Lemma \ref{compact} one obtains, for any $\mu$-inner function $\theta \in I_{\mu}$,
\[
T^*_{\theta} X T_{\theta} = X + {\rm WOT}-\lim_{k \rightarrow \infty} T^*_{\theta_k}(T^*_{\theta} S T_{\theta} - S)T_{\theta_k}
= X + |\alpha|^2(T^*_{\theta} S T_{\theta} - S),
\]
or equivalently,
\[
T^*_{\theta}(X - |\alpha|^2 S) T_{\theta} = X - |\alpha|^2 S.
\]
Hence there is a symbol $Y \in (U)'$ with
\[
T_Y = \frac{1}{1 - |\alpha|^2}(X - |\alpha|^2 S).
\]
Thus $S$ is a Schatten-$p$-class perturbation of an abstract $T$-Toeplitz operator if and only if $S - X \in \mathcal S_p(\Hcal)$. 
Now let us suppose that $1 \leq p < \infty$. Since
\[
S - X = {\rm WOT}-\lim_{k \rightarrow \infty} (S - T^*_{\theta_k} S T_{\theta_k})
\]
and since the Schatten-$p$-norm is lower semicontinuous with respect to the weak operator topology, to conclude that
$S - X \in \mathcal S_p(\Hcal)$, it suffices to choose the sequence $(\theta_k)_{k \in \mathbb N}$ in such a way that 
\[
\sup_{k \in \mathbb N} \| S - T^*_{\theta_k} S T_{\theta_k} \|_p < \infty.
\]
Our results on products of inner functions and the following extension of the $\epsilon$-$\delta$-criterion from \cite{xia-characterization}
will allow us to complete the proof along these lines.

\begin{prop} \label{eps-delta}
Let $1 \leq p \leq \infty$ and $w \in D$ be arbitrary. Suppose that $S \in B(\Hcal)$ is an operator such that
\[
T^*_{\theta} S T_{\theta} - S \in \mathcal S_p(\Hcal)
\]
for all $\theta \in I_D$. Then, for every real number $\epsilon > 0$, there is a real number $\delta > 0$ such that
\[
\| T^*_{\theta} S T_{\theta} - S \|_p < \epsilon
\]
for each function $\theta \in I_D$ with $1 - |\theta(w)| < \delta$.
\end{prop}

\proof Let us assume that the assertion is not true. Then there are a positive real number $\epsilon > 0$
and a sequence $(\theta_k)_{k \in \mathbb N}$ in $I_D$ with $1 - |\theta_k(w)| < 2^{-k}$ and 
\[
\| T^*_{\theta_k} S T_{\theta_k} - S \|_p \geq  \epsilon
\]
for all $k \in \mathbb N$. Multiplying the functions $\theta_k$ with suitable complex numbers of modulus one, we
can achieve that $\theta_k(w) > 0$ for all $k$.

Using part (i) and (iii) of Corollary \ref{xia1to3} and Lemma \ref{SOT}, and arguing otherwise exactly as in the
proof of Lemma 5 from \cite{xia-characterization}, one can show that there are a strictly increasing sequence
$(k(j))_{j \geq 1}$ of natural numbers and a sequence $(F_j)_{j \geq 1}$ of pairwise orthogonal finite-rank orthogonal
projections such that the partial products
\[
v_j = \prod_{\nu = 1}^j \theta_{k(\nu)} \quad \quad (j \geq 1)
\]
satisfy the inequalities
\begin{enumerate}[(i)]
\item $\| F_{j+1} W_j F_{j+1} \|_p \geq \epsilon /2$,
\item $\| (1 - F_{j+1}) W_j \|_p \leq 2^{-j}$,
\item $\| W_j (1 - F_{j+1}) \|_p \leq 2^{-j}$
\end{enumerate}
for all $j \geq 1$, where by definition
\[
W_j = T^*_{v_j} (T^*_{\theta_{k(j+1)}} S T_{\theta_{k(j+1)}} - S) T_{v_j} \quad (j \geq 1).
\]
Then, for each $j \geq 1$, the operator
\[
K_j = W_j - F_{j+1} W_j F_{j+1} \in \mathcal S_p(\Hcal)
\]
satisfies the estimate
\[
\| K_j \|_p \leq \| (1 - F_{j+1}) W_j \|_p + \| F_{j+1} W_j (1 - F_{j+1}) \|_p \leq 2^{-j+1}.
\]
Hence $K = \sum_{j=1}^{\infty} K_j \in \mathcal S_p(\Hcal)$ is a well-defined operator. Since the
orthogonal projections $F_j$ $(j \geq 1)$ are pairwise orthogonal, the series
\[
B = \sum_{j=1}^{\infty} B_j \quad {\rm with} \quad B_j =  F_{j+1} W_j F_{j+1} \in B(\Hcal)
\]
converges in the strong operator topology. By construction, $B \notin S_p(\Hcal)$. Indeed, otherwise
Lemma \ref{biSOT} would lead to the contradiction
\[
\frac{\epsilon}{2} \leq \| F_{j+1} W_j F_{j+1} \|_p = \| F_{j+1} B F_{j+1} \|_p \leq \| F_{j+1} B \|_p \stackrel{k}{\rightarrow} 0.
\]
Since $ \sum_{j=1}^{\infty} |1 - \theta_{k(j)}(w) | < \infty$, it follows from Corollary \ref{prod2} that the sequence $(v_j)_{j \geq 1}$
is weak$^*$ convergent in $H^{\infty}(D)$ to the inner function 
\[
v = \prod_{j=1}^{\infty} \theta_{k(j)} \in I_D
\]
and in combination with Lemma \ref{biSOT} we obtain that
\[
{\rm SOT}-\lim_{j\rightarrow \infty} T_{v_j} = T_v \quad {\rm and} \quad {\rm SOT}-\lim_{j\rightarrow \infty} T^*_{v_j} = T^*_v
\]
in $B(\Hcal)$. By construction we have
\[
K + B = {\rm SOT}-\sum_{j=1}^{\infty} W_j = {\rm SOT}-\sum_{j=1}^{\infty} T^*_{v_j} (T^*_{\theta_{k(j+1)}} S T_{\theta_{k(j+1)}} - S) T_{v_j}
\]
\[
= {\rm SOT}-\sum_{j=1}^{\infty} (T^*_{v_{j+1}} S T_{v_{j+1}} - T^*_{v_j} S T_{v_j})
\]
\[
= (T^*_v S T_v - S) - (T^*_{v_1} S T_{v_1} - S) \in \mathcal S_p(\Hcal).
\]
Since also $K \in \mathcal S_p(\Hcal)$, we obtain that $B \in \mathcal S_p(\Hcal)$. This contradiction completes the proof.
\proofend

After these preparations we can prove our main result on Schatten class perturbations of Toeplitz operators.
As before, let $T \in B(\Hcal)^n$ be a regular $A$-isometry with respect to the domain algebra $A = A(D)$ of a
bounded domain $D \subset \mathbb C^n$ such that the associated scalar spectral measure $\mu \in M^+_1(\partial_A)$
is a faithful Henkin measure.

\begin{thm} \label{perturbation}
Let $1 \leq p < \infty$ and let $S \in B(\Hcal)$. Then
\[
T^*_{\theta} S T_{\theta} - S \in \mathcal S_p(\Hcal) \; {\rm for  \; all} \; \theta \in I_D
\]
if and only if $S \in \Tcal(T) + \mathcal S_p(\Hcal)$.
\end{thm}

\proof Let $S \in B(\Hcal)$ be an operator such that $T^*_{\theta} S T_{\theta} - S \in \mathcal S_p(\Hcal)$ for all $\theta \in I_D$.
By the remarks preceding Proposition \ref{eps-delta} it suffices to show that $S \in \Tcal(T) + \mathcal S_p(\Hcal)$. Fix an arbitrary
point $w \in D$. By Proposition \ref{eps-delta} there is a number $\delta \in (0,1)$ such that
\[
\| T^*_{\theta} S T_{\theta} - S \|_p < 1
\]
for all $\theta \in I_D$ with $1 - |\theta(w)| < \delta$. Define $\alpha = 1 - \delta/2$. Since $(A(D)|\partial_A,\partial_A,\mu)$ is regular in
the sense of Aleksandrov and since $\mu$ is continuous, there is a sequence $(\theta_k)_{k \in \mathbb N}$ in $I_D$ with
$w^*-\lim_{k \rightarrow \infty} \theta^*_k = \alpha$. Since
\[
r_{\mu}: H^{\infty}(D) \rightarrow L^{\infty}(\mu)
\]
has weak$^*$ closed range and since this map induces a weak$^*$ homeomorphisms onto its range $\ran r_{\mu} \supset H_A^{\infty}(\mu)$, it follows
that $w^*-\lim_{k \rightarrow \infty} \theta_k = \alpha$ in $H^{\infty}(D)$. By passing to a subsequence we can achieve that
$1 - |\theta_k(w)| < \delta$ for all $k$ and that at the same time the limit
\[
X = {\rm WOT}-\lim_{k \rightarrow \infty} T^*_{\theta_k} S T_{\theta_k} \in B(\Hcal)
\]
exists. Then the remarks preceding Proposition \ref{eps-delta} show that $S \in \Tcal(T) + \mathcal S_p(\Hcal)$.
\proofend

It follows from the results in Section 5 of \cite{subnormal} that the canonical probability measure $\sigma$ on the Shilov boundary of a
smooth strictly pseudoconvex or bounded symmetric and circled domain $D$ in $\mathbb C^n$ is a faithful Henkin measure. Since $\sigma$ is a
scalar spectral measure of the minimal normal extension $M_z \in B(L^2(\sigma))^n$ of the regular $A$-isometry $T_z \in B(H^2(\sigma))^n$
and since $W^*(M_z) \subset B(L^2(\sigma))$ is a maximal abelian $W^*$-algebra, Theorem \ref{xia} follows as an immediate application of 
Theorem \ref{perturbation}.

M. Didas, J. Eschmeier, D. Schillo\\
Fachrichtung Mathematik\\
Universit\"at des Saarlandes\\
Postfach 151150\\
D-66041 Saarbr\"ucken\\
Germany

\texttt{didas@math.uni-sb.de, eschmei@math.uni-sb.de, schillo@math.uni-sb.de}

\end{document}